\documentclass[11pt, a4paper]{article}
\usepackage{amsmath,amssymb,latexsym,graphics,epsfig}
\usepackage{hyperref}
\usepackage{color}
\usepackage{amsthm}

\newtheorem{theo}{Theorem}[section]
\newtheorem{lema}[theo]{Lemma}

\newtheorem{coro}[theo]{Corollary}

\def\proof{{\boldmath  $Proof.$}\hskip 0.3truecm}

\def\endproof{\mbox{ \quad $\Box$}}

\def\vecalpha{\mbox{\boldmath $\alpha$}}
\def\vecbeta{\mbox{\boldmath $\beta$}}
\def\vecgamma{\mbox{\boldmath $\gamma$}}
\def\vecdelta{\mbox{\boldmath $\delta$}}
\def\vecsigma{\mbox{\boldmath $\sigma$}}

\def\vecmu{\mbox{\boldmath $\mu$}}

\begin{document}
\title{A Note on the Voting Problem
}

\author{M.A. Fiol
\\ \\
{\small Universitat Polit\`ecnica de Catalunya, BarcelonaTech}\\
{\small Departament de Matem\`atica Aplicada IV} \\
{\small Barcelona, Catalonia} \\
{\small (e-mail: {\tt fiol@ma4.upc.edu})}}

\date{}

\maketitle

\begin{abstract}
Let $v(n)$ be the minimum number of voters with transitive preferences which are needed to generate any strong preference pattern (ties not allowed) on $n$ candidates. Let $k=\lfloor \log_2 n\rfloor$. We show that $v(n)\le n-k$ if $n$ and $k$ have different parity, and $v(n)\le n-k+1$ otherwise.
\end{abstract}

\section{Introduction}
Let us consider a set of $n$ candidates or options $A=\{a,b,c,\ldots\}$ which are ordered by order of preference by each individual of a set $U$ of voters. Thus, each $\vecalpha\in U$ can be identified with a permutation $\vecalpha=x_1x_2\cdots x_n$ of the elements of $A$, where $x_i$ is preferred over $x_j$ (denoted $x_i\rightarrow x_j$) if and only if $i<j$. The set of voters determine what is called a {\em preference pattern} which summarizes the majority opinion about each pair of options.

In this note only {\em strong} preference patterns are considered, that is, it is assumed that there are no ties. So, each preference pattern on $n$ options is fully represented by a tournament $T_n$ on $n$ vertices where the arc $(a,b)$ means $a\rightarrow b$, that is, $a$ is preferred over $b$ by a majority of voters.
Conversely, given any pattern $T_n$ we may be interested in finding a minimum set of voters, denoted $U(T_n)$, which generates $T_n$. Let $v(T_n)=|U(T_n)|$ and let $v(n)=\max\{v(T_n)\}$ computed over all tournaments with $n$ vertices. In \cite{mcg53} McGarvey showed that $v(n)$ is well defined, that is, for any $T_n$ there always exist a set $U(T_n)$ and $v(n)\le 2{n \choose 2}$. Sterns \cite{s59} showed that $v(n)\le n+2$ if $n$ is even and $v(n)\le n+1$ if $n$ is odd. Finally, Erd\"{o}s and Moser \cite{em64} were able to prove that $v(n)$ is of the order $O(n/\log_2 n)$. In fact all the above results were given for preference patterns which are not necessarily strong (in this case a tie between $a$ and $b$ can be represented either by an absence of arcs between $a$ and $b$ or by an edge $\{a,b\}$). It is worth noting that, contrarily to the method of Erd\"{o}s and Moser, the approaches of McGarvey and Sterns give explicit constructions of a set of voters which generate any desired pattern. In the case of strong patterns we improve the results of the latter authors by giving and inductive method to obtain a suitable set of voters.

\section{Strong preference patterns}
Let us begin with a very simple but useful result, which is a direct consequence of the fact that in our preference patterns there are no ties.

\begin{lema}
\label{basiclema}
Let $v(n)$ be defined as above. Then, $v(n)$ is odd.
\end{lema}

\proof
By contradiction, suppose that, for a given strong pattern $T_n$, $v(n)$ is even. Then, for any two options $a,b$ we have that either $a\rightarrow b$ or $b\rightarrow a$ with at least two votes of difference. Consequently, the removing of a voter does not change the preference pattern.
\endproof

Notice that, from this lemma, Sterns' result particularized for strong patterns are $v(n)\le n+1$ for $n$ even and $v(n)\le n$ for $n$ odd.

Our results are based on the following theorem.

\begin{theo}
\label{theo1}
Let $T_{n+2}$ be a strong pattern containing two options, say $a$ and $b$. Let $T_{n}=T_{n+2}\setminus\{a,b\}$. Then,
$v(T_{n+2})\le v(T_n)+2$.
\end{theo}

\proof
Let $U(T_n)=\{\vecalpha_1,\vecalpha_2,\ldots,\vecalpha_r\}$ be a minimum set of $r=v(T_n)$ voters generating $T_n$. By Lemma \ref{basiclema}, $r$ is odd. Besides, suppose without loss of generality that $a\rightarrow b$, and consider the sets $A_1=\{x\neq a\,|\, x\rightarrow b\}$ and $A_2=\{x\neq b\,|\, a\rightarrow x\}$. Assuming $A_1\cap A_2\neq \emptyset,A_1,A_2$ (any other case follows trivially from this one), we can write $A_1=\{y_1,y_2,\ldots,y_s,\ldots,y_t\}$ and $A_2=\{y_s,y_{s+1},\ldots,y_t,\ldots,y_m\}$, $1<s\le t<m$. Now, let us define the sequences $\vecgamma=y_1y_2\cdots y_{s-1}$, $\vecdelta=y_sy_{s+1}\cdots y_t$, $\vecsigma=y_{t+1}y_{t+2}\cdots y_m$ and
$\vecmu=y_{m+1}y_{m+2}\cdots y_n$, and consider the following set of $r+2$ voters:
\begin{eqnarray*}
\vecbeta_i&=& b\vecalpha_i a,\qquad 1\le i\le (r+1)/2, \\
\vecbeta_j&=& a\vecalpha_j b,\qquad (r+3)/2\le j\le r, \\
\vecbeta_{r+1} &=& \vecgamma a \vecdelta b\vecsigma \vecmu, \\
\vecbeta_{r+2} &=& \overline{\vecmu} a\overline{\vecsigma} \overline{\vecdelta} \overline{\vecgamma}b,
\end{eqnarray*}
where $\overline{\vecgamma}=y_{s-1}\cdots y_2 y_1$, $\overline{\vecdelta}=y_t\cdots y_{s+1}y_s$, etc. Now it is routine to verify that these voters generate the pattern $T_{n+2}$ and, hence, $v(T_{n+2})\le r+2=v(T_n)+2$.
\endproof

A tournament or strong preference pattern $T$ is called {\em transitive} if $a\rightarrow b$ and $b\rightarrow c$ implies $a\rightarrow c$. In this case it is clear that $v(T)=1$. The proof of the following result can be found in \cite{em64}.

\begin{theo}[\cite{em64}]
\label{theo2}
Let $f(n)$ be the maximum number such that every tournament on $n$ vertices has a transitive subtournament on $f(n)$ vertices. Then,
$$
\lfloor \log_2 n\rfloor + 1\le f(n)\le 2\lfloor \log_2 n\rfloor+1.
$$
\end{theo}

The proof of the lower bound, due to Sterns, gives a very simple algorithm to find a subtournament which attains such a bound, see again \cite{em64}.

From Theorems \ref{theo1} and \ref{theo2} we get the following corollary.

\begin{coro}
Given $n\ge 2$, set $k=\lfloor\log_2 n\rfloor$. Then $v(n)\le n-k$ if $n$ and $k$ have different parity, and $v(n)\le n-k+1$ otherwise.
\end{coro}

\proof
Let $T_n$ be any tournament on $n$ vertices. First, use Theorem \ref{theo2} to find a transitive subtournament $T$ on $k+1$ vertices. If $n$ and $k$ have different parity, then $n-k-1$ is even. So, starting from $T$, we can apply Theorem \ref{theo1} repeatedly, $(n-k-1)/2$ times, to obtain a set of $n-k$ voters which generates $T_n$. Otherwise, we consider a subtournament of $T$ on $k$ vertices and proceed as above with the remaining $n-k$ vertices. \endproof


\end{document}